\documentclass[10pt,twoside]{amsart}
\usepackage{amssymb,epsfig,amsmath,latexsym,amscd}
\usepackage{tikz}
\usepackage{enumitem}
\usepackage{comment}
\nonstopmode

\numberwithin{equation}{section}
\font\tengothic=eufm10 scaled\magstep 1
\font\sevengothic=eufm7 scaled\magstep 1
\newfam\gothicfam
\textfont\gothicfam=\tengothic
\scriptfont\gothicfam=\sevengothic

\DeclareMathOperator{\pnt}{\raise 0.5mm \hbox{\large\bf.}}

\newtheorem{theorem}{Theorem}%[section]

\theoremstyle{definition}
\newtheorem{definition}[theorem]{Definition} % \theoremstyle{remark}

\newtheorem{example}[theorem]{Example}

\author[G. Favacchio]{Giuseppe Favacchio}
\address{Dipartimento di Ingegneria, Viale delle Scienze - 90128 - Palermo, Italy}
\email{giuseppe.favacchio@unipa.it}

\author[E. Guardo]{Elena Guardo}
\address{Dipartimento di Matematica e Informatica, Viale A. Doria, 6 - 95100 - Catania, Italy}
\email{guardo@dmi.unict.it}

\author[J.~Migliore]{Juan Migliore}
\address{Department of Mathematics, University of Notre Dame, Notre Dame, IN 46556}
\email{migliore.1@nd.edu}

\title[Correction to: Proc. Amer. Math. Soc. 146 (2018), 2811-2825]{Correction to: On the arithmetically Cohen-Macaulay property for sets of points in multiprojective spaces. Proc. Amer. Math. Soc. 146 (2018), 2811-2825}

\begin{document}

\thanks{Revised Version:  August 20, 2024}
\begin{abstract}
We correct a mistake in the cited paper. It introduced a combinatorial property, the $(\star_n)$-property, for a finite set of points $X$ in $(\mathbb P^1)^n$ and claimed that this property holds if and only if $X$ is ACM. In fact $X$ being ACM is a sufficient condition for the $(\star_n)$-property, but we only prove that it is necessary when $n=3$, and we give a counterexample when $n=4$.
\end{abstract}
\maketitle

    The purpose of this corrigendum is to correct a mistake in \cite[Theorem 3.16]{FGM}, which claimed that a set of points $X$ in $(\mathbb P^1)^n$ is arithmetically Cohen-Macaulay (ACM) if and only if a certain property called $(\star_n)$ holds. All other results remain true as they were stated. We thank Gunnar Fløystad for pointing out a gap in our proof and providing a counterexample.

In \cite[Theorem 3.16]{FGM},  part ($\sigma_1$), first bullet we claimed that $\hat{Y}_1 \cap Y_2$ contains no subset of type (i). 
In particular, in the third paragraph we wrote:

{\em

Then by applying \cite[Lemma 3.15]{FGM} to the points $P_{(1,v_2,
\ldots, v_n)}$ and $P_{\underline u_{s-1}}$, we get a
contradiction by forcing a point of $\hat{Y}_1 \cap Y_2$ to lie in
the complete intersection.
}

This fact doesn't lead necessarily to a contradiction since the mentioned point of $\hat{Y}_1 \cap Y_2$ forced to lie in
the complete intersection could be $P_{\underline{v}}.$ 

However, the gap in the proof of \cite[Theorem 3.16]{FGM} cannot be fixed. Indeed, G. Fløystad provided us the following example.
\begin{example}
Let $R=K[x_1,y_1,x_2,y_2, \cdots, x_4,y_4]$ be the coordinate ring of $(\mathbb P^1)^4$. 
Let $A=[0:1], B=[1:0]\in \mathbb P^1$ and consider the two sets of four points in $(\mathbb P^1)^4$ $$Y_1= \{(A,A,A,A), (B,A,A,A), (B,B,A,A), (B,B,B,A)\}$$ and $$Y_2=\{(B,B,B,B), (A,B,B,B), (A,A,B,B), (A,A,A,B)\}.$$
Set $X=Y_1\cup Y_2$.  
Then $\hat{Y}_1 \cap Y_2= \{(B,B,B,B), (A,A,A,B)\}$
which is not ACM. 

Moreover, one can check that $X$ has the $(\star_4)$ property but
 $R/I_X$ is not ACM. In particular 
$$I_X=(x_1y_1, x_2y_2, x_3y_3, x_4y_4, x_1x_3y_2, x_1x_4y_2, x_1x_4y_3, x_2x_4y_3, x_2y_1y_3, x_2y_1y_4, x_3y_1y_4, x_3y_2y_4)$$
and the Betti table of $R/I_X$ is 
\[
\begin{matrix}
 & 0 & 1 & 2 & 3 & 4 & 5 & 6\\
0: & 1 & . & . & . & . & . & .\\
1: & . & 4 & . & . & . & . & .\\
2: & . & 8 & 38 & 48 & 28 & 8 & 1
\end{matrix}
\]
\end{example}

Thus we need to weaken the statement of Theorem 3.16, and we replace it with the following result.

\begin{theorem}\label{mainCharacterization2}
    Let $X \subset (\mathbb P^1)^n $ be a finite set. If  $X$ is ACM then $X$ has the $(\star_n)$-property.
\end{theorem}

The original proof was correct for this part. For the reverse implication, Theorem 3.16 is still true for points in $(\mathbb P^1)^3.$ For the convenience of the reader, we write here the original proof adapted for the case $(\mathbb P^1)^3$. However, only the item $(\sigma_1)$ needed to be fixed; we believe that the rest is correct. 

\begin{theorem}\label{mainCharacterization}
   Let $X \subset \mathbb P^1 \times \mathbb P^1 \times \mathbb P^1$ be a finite set. Then $X$ has the $(\star_3)$-property if and only if $X$ is ACM.
\end{theorem}
\begin{proof}
As noted above, 
if $X$ is ACM, it was shown in \cite[Corollary 3.11]{FGM}  that even more generally in $(\mathbb P^1)^n$, the $(\star_n)$ property holds.
   For the converse,  we assume $n=3$ and 
   we proceed by induction on $t$, the number of level sets with respect to some projection (say $\pi_1$, { i.e., the projection omitting the first component, see \cite[Definition 2.1]{FGM}}). If $t$ is equal to 1  the result follows from \cite[Corollary 2.9]{FGM}. Let $X = X_1 \cup X_2\cup \cdots \cup X_t$   { be the natural stratification of the points of $X$ as the union of $1$-level sets (see \cite[Definition 2.5]{FGM}).
   Let $Y_1 = X_1$ be one of these 1-level sets (this could be any of them up to re-indexing)} and let $Y_2 = X_2 \cup \dots \cup X_t$ be the union of the remaining level sets. To shorten the notation, in the sequel of the proof, we denote by $A_i\in R_{(1,0,0)}$, $B_j\in R_{(0,1,0)}$ and $C_k\in R_{(0,0,1)}$ the linear forms defining the multihomogeneous hyperplanes. Moreover we denote by $P_{ijk}$ the point whose ideal is generated by $(A_i,B_j,C_k)$ and by $L_{jk}$ the line of type $(0,1,1)$ whose ideal is generated by $(B_j,C_k),$ for some $A_i\in R_{(1,0,0)}$, $B_j\in R_{(0,1,0)}$ and $C_k\in R_{(0,0,1)}.$ We assume $A_1\in R_{(1,0,0)}$ is the linear form defining the hyperplane containing $Y_1$.
   We denote by $\hat{Y}_1$ the set of lines $L_{jk}$ passing through   points of $Y_1$ (one for each line),   i.e. $\hat{Y}_1 = \pi_1^{-1}(\pi_1(X_1))$. (Viewed in $\mathbb P^5$, $\hat{Y}_1$ is a union of codimension 2 linear spaces.)

   We know by induction that both $Y_1$ and $Y_2$ are ACM. In particular, this means that $\hat{Y}_1$ is also ACM. Hence we have an equality of saturated ideals $I_{Y_1}=(A_1)+ I_{\hat{Y}_1}$. Then it follows from the following exact sequence
\[
0 \rightarrow I_{{Y}_1} \cap I_{{Y}_2} \rightarrow I_{{Y}_1} \oplus I_{{Y}_2} \rightarrow (A_1)+ I_{\hat{Y}_1} + I_{{Y}_2} \rightarrow 0
\]
that it is enough to show that $I_{\hat{Y}_1} + I_{{Y}_2}$ is an ACM ideal (hence clearly of height 3) and $A_1$ is a regular form in $R/(I_{\hat{Y}_1} + I_{{Y}_2})$. We proceed by steps.

\begin{itemize}[leftmargin=0.8cm]

\item[$(\sigma_1)$]
{\em We show that $\hat{Y}_1\cap{{Y}_2}$ is an ACM set of points. }

\smallskip

By the inductive hypothesis, it suffices to show that it has the $(\star_3)$-property, i.e. $\hat{Y}_1\cap{{Y}_2}$ does not contain any of the configurations (a), (b), (c) in \cite[Remark 3.10]{FGM}.

\smallskip

\begin{itemize}

\item[$\bullet$] {\em $\hat{Y}_1 \cap Y_2$ contains no subset of type (a)}. First consider another level set with respect to $\pi_1$; without loss of generality assume it is the second one. If $P_{211}, P_{222}\in {\hat{Y}_1\cap {Y}_2},$ then since $X$ does not have a subset of type (a), it follows that either $P_{212}\in X$ or $P_{221}\in X.$  Assume without loss of generality that $P_{212}\in X.$ If $P_{112}\not\in X,$ then since $P_{111}, P_{122}\in Y_1$, which is ACM, it follows that  $P_{121}\in X$. But $X$ does not contain a subset of type (c), so $P_{221}\in X,$ and therefore $P_{221}\in \hat{Y}_1\cap{{Y}_2}.$

\smallskip

Moreover suppose that a level set of  $\hat{Y}_1\cap{{Y}_2}$ with respect to another direction contains another subset of type (a);  say $\{P_{211}, P_{321}\}\in \hat{Y}_1\cap{{Y}_2}$ and $\{P_{221}, P_{311}\}\notin \hat{Y}_1\cap{{Y}_2}$. Then also $X$ must contain the same configuration, giving a contradiction.

\smallskip

\item[$\bullet$] {\em $\hat{Y}_1 \cap Y_2$ contains no subset of type $(b)$};
Indeed, suppose  that  ${\hat{Y}_1\cap {Y}_2}$  contains a subset of type $(b)$; say this subset is \begin{equation}\label{eq.1}
    P_{211}\in \hat{Y}_1\cap {Y}_2,\ \ \ P_{322}\in \hat{Y}_1\cap {Y}_2
\end{equation}   and
\begin{equation}\label{eq.2}
    P_{212}, P_{221}, P_{321}, P_{312}\notin \hat{Y}_1\cap {Y}_2.
\end{equation}
Moreover, from \eqref{eq.1}, we get

\begin{equation}\label{eq.3}P_{111}\in X\ \ \ \text{and}\ \ \  P_{122}\in X.\end{equation} 
Since $X$ has the $(\star_3)$-property, from \eqref{eq.3} either $P_{121},$ or $P_{112}\in Y_1$; and there is a path (obtained by
changing one coordinate at a time) in $Y_2$ joining $P_{211}$ and $P_{322}$.
We note that if both $P_{121},P_{112}\in X$ then we get a contradiction since then a path joining $P_{211}$ and $P_{322}$ in $X$ is also contained in $\hat{Y}_1\cap {Y}_2$. 
So, we assume (without loss of generality) $$P_{121}\in X\ \ \ \text{and}\ \ \ P_{112}\notin X.$$ The possible paths in $X$ connecting $P_{211}$ and $P_{322}$ are listed below. We show that each case leads to a contradiction.   
\begin{itemize}
    \item[1.] If $P_{221},$ $P_{321}\in X$, since $P_{121}\in Y_1$, then these two points also are in $\hat{Y}_1\cap {Y}_2$ contradicting the condition in \eqref{eq.2};
    \item[2.] If $P_{221},$ $P_{222}\in X$, since $P_{121},P_{122}\in Y_1$ then these two points also are in $\hat{Y}_1\cap {Y}_2$ contradicting the condition in \eqref{eq.2};
    \item[3.] If $P_{311},$ $P_{321}\in X$ since $P_{111},P_{121}\in Y_1$ then these two points also are in $\hat{Y}_1\cap {Y}_2$ contradicting the condition in \eqref{eq.2};
    \item[4.] If $P_{311},$ $P_{312}\in X$. The set  $\{P_{111}, P_{121}, P_{122},P_{311}, P_{312}, P_{322}\} \subseteq X$ is of type $(c)$. Since $X$ has the $(\star_3)$-property and $P_{112}\notin X$ we get $P_{321}\in X.$ Then we are in the case of item 3.   
    \item[5.] If $P_{212},$ $P_{222}\in X$.  The set  $\{P_{111}, P_{121}, P_{122},P_{211}, P_{212}, P_{222}\} \subseteq X$ is of type $(c)$. Since $X$ has the $(\star_3)$-property and $P_{112}\notin X$ we get $P_{221}\in X.$ Then we are in the case of item 2.    
    \item[6.] If $P_{212},$ $P_{312}\in X$. The set  $\{P_{122}, P_{212} \} \subseteq X$ is of type $(a)$. Since $X$ has the $(\star_3)$-property and $P_{112}\notin X$ we get $P_{222}\in X.$ Then we are in the case of item 5.    
\end{itemize}

\smallskip

\item[$\bullet$] {\em $\hat{Y}_1 \cap Y_2$ contains no subset of type (c)}. Indeed, if it did then this subset is contained in $X$, contradicting the $(\star_3)$-property of $X$.

\end{itemize}

\medskip

\item[$(\sigma_2)$] {\em We make a technical observation concerning the ``outlier" points.}

\smallskip

We denote by $Y'_1$ the set of points $P_{1jk}\in Y_1$ such that the line $L_{jk}$ has empty intersection with ${Y}_2,$ and we denote by $Y_2':=Y_2\setminus (\hat{Y}_1\cap {Y}_2).$  Let $F\in I_{\hat{Y}_1\cap {Y}_2}$ be a product of linear forms of type $A_i$, $B_j$ and $C_k$. Taking the ideal of the empty set to be $R$, {\em we claim that }
\[
F\in (I_{Y_2'}) \cup (I_{Y_1'}).
\]
We assume that both $Y_1'$ and $Y_2'$ are non-empty; otherwise the statement is trivial.
   Assume by contradiction that $F\notin (I_{Y_2'}) \cup (I_{Y_1'}).$ Then there exist $P_{111}\in Y_1'$ (since $X$ has the $(\star_3)$-property and by the definition of these two sets, this implies $P_{212},P_{221}\notin Y_2'$) and $P_{222}\in Y_2',$
   such that $F\notin(A_1,B_1,C_1)$ and $F\notin(A_2,B_2,C_2)$.  Now, $P_{111}\in Y_1'$ implies $P_{211}\notin Y_2$;  moreover, since $P_{222}\in Y_2'$ we have $P_{122}\notin Y_1$. (Here we use the fact that since $X$ has the $(\star_3)$-property, this excludes either $B_1=B_2$ or $C_1=C_2$.) But $X$ does not contain subsets of type (a) or (b) so, without loss of generality, we can assume $P_{121},P_{221}\in X$ i.e. $P_{221}\in {\hat{Y}_1\cap {Y}_2}.$ Since $F\in (A_2,B_2,C_1)$ and it is a product of linear forms, we get either $F\in(A_2,B_2,C_2)$ or $F\in(A_1,B_1,C_1)$, which contradicts the assumption.

\medskip

    \medskip

\item[$(\sigma_3)$] {\em We show that $I_{\hat{Y}_1\cap {Y}_2}\subseteq I_{\hat{Y}_1}+I_{Y_2}$.}

\smallskip

From $(\sigma_1)$ we know that $\hat{Y}_1\cap {Y}_2$ is ACM, so  $I_{\hat{Y}_1\cap {Y}_2}$ is minimally generated by products of linear forms (from \cite[Corollary 3.4]{FGM}).  Let $F\in I_{\hat{Y}_1\cap {Y}_2}$ be such a generator. From the minimality of $F$ we note that $F\notin (A_1)$.
From $(\sigma_2)$ we have $F\in I_{Y_1'}\cup I_{Y_2'}$. Assume first that $F\in I_{Y_2'}$. Then,  from the definition of $Y_2'$, trivially we have $F\in I_{Y_2}\subseteq I_{\hat{Y}_1}+ I_{Y_2}.$

Assume now that $F\notin I_{Y_2'}$ i.e. there exists a point, say $P_{222}\in {Y_2'}$, such that $F\notin(A_2,B_2,C_2)$. We collect the relevant facts:

\medskip

\begin{itemize}
\item[$(f_1)$] $P_{122} \notin Y_1$ by definition of $Y_2'$; \vspace{.07in}

\item[$(f_2)$] $F \in I_{Y_1'}$ by $(\sigma_2)$; \vspace{.07in}

\item[$(f_3)$] $F \in I_{\hat{Y}_1 \cap Y_2}$ but $F \notin (A_2,B_2,C_2)$ and $F \notin (A_1)$ by assumption.
\end{itemize}

\medskip

\noindent We want to show that $F\in I_{\hat{Y}_1}.$ Choose any point  $P \in Y_1$. Suppose first that $P = P_{121}$. Then in order to avoid $X$ having a subset of type (a), from $(f_1)$ we must have $P_{221} \in X$. This means that $P_{121} \notin Y_1'$, so $P_{221}$  lies in $\hat{Y}_1 \cap Y_2$.  thus, We have $F(P_{221}) = 0$  but $F(P_{222}) \neq 0$. It follows that $F \in (B_2,C_1) \subseteq I_{\hat Y_1}$.

The case $P = P_{112}$ is entirely analogous.  Finally, suppose that $P$ is any other point of $Y_1$; without loss of generality say it is $P_{111}$.
We consider two cases.

\smallskip

\begin{itemize}
\item[$\bullet$] If $P_{111}\in Y_1'$ then from $(f_2)$ we have $F\in (A_1,B_1,C_1)$ and, since $F \notin (A_1)$,  $F\in (B_1,C_1)\subseteq I_{\hat Y_1}$.

\smallskip

\item[$\bullet$] If $P_{111}\in Y_1\setminus Y_1'$ then since $X$ has the $(\star_3)$-property, to avoid a subset of type (a) or (b),  at least one among $P_{211}$, $P_{212}$ and $P_{221}$ belongs to $\hat{Y}_1\cap {Y}_2.$  Therefore, since $F\in I_{\hat{Y}_1\cap {Y}_2}$ and $F\notin(A_2,B_2,C_2)$ we get $F\in (B_1,C_1)\subseteq I_{\hat Y_1}.$

\end{itemize}

\noindent This shows that $F \in I_{\hat{Y}_1}$ as desired, and concludes the proof of $(\sigma_3)$.

\end{itemize}

\medskip

To complete the proof of our theorem, note that  on the other hand we always  have $I_{\hat{Y}_1\cap {Y}_2} = \sqrt{I_{\hat{Y}_1} + I_{Y_2}}\supseteq I_{\hat{Y}_1}+I_{Y_2}.$ Thus, $I_{\hat{Y}_1}+I_{Y_2}$ is the ideal of  an ACM set of reduced points $({\hat{Y}_1\cap {Y}_2}$) in $\mathbb{P}^1\times\mathbb{P}^1\times\mathbb{P}^1$, as desired. Moreover, this implies that $A_1$ is a regular form in $R/(I_{\hat{Y}_1}+I_{Y_2})$ since no point of $\hat{Y}_1\cap {Y}_2$  belongs to the plane defined by $A_1.$
\end{proof}

\end{document}